\documentclass[12pt,reqno]{amsart}

\newcommand{\C}{\mathbb{C}}\newcommand{\Q}{\mathbb{Q}}
\newcommand{\N}{\mathbb{N}}\newcommand{\R}{\mathbb{R}}
\newcommand{\Z}{\mathbb{Z}}

\newcommand{\G}{{\mathbb G}}

\newcommand{\ta}{\theta}\newcommand{\disc}{\operatorname{disc}}
\newcommand{\eps}{\varepsilon}
\newcommand{\al}{\alpha}\newcommand{\om}{\omega}
\newcommand{\be}{\begin{equation}}
\newcommand{\ee}{\end{equation}}
\newcommand{\bes}{\begin{equation*}}
\newcommand{\ees}{\end{equation*}}

\usepackage{amsfonts}
\usepackage{amssymb}
\usepackage{enumerate}
\usepackage{comment}
\usepackage{epsfig}

\def\tareesidedbox#1{\setbox0=\hbox{$#1$}\dimen0=\wd0 \advance\dimen0 by3pt\rlap{\hbox{\vrule height8pt width.4pt
 depth2pt \kern-.4pt\vrule height8.4pt width\dimen0 depth-8pt\kern-.4pt \vrule height8pt width.4pt depth2pt}}
 \relax \hbox to\dimen0{\hss$#1$\hss}}
\def\ho#1{\tareesidedbox#1}
\def\eps{\varepsilon}

\title[   Mahler measure of algebraic numbers]{
The Mahler measure of algebraic numbers:\\
a survey}

\author{Chris Smyth }

\begin{document}

\begin{abstract}
A survey of results for Mahler measure of algebraic numbers, and
  one-variable polynomials with integer coefficients is presented. Related results on the maximum modulus of the
  conjugates (`house') of an algebraic integer are also
discussed. Some generalisations are given too, though not to Mahler measure of polynomials in more than one
variable.
\end{abstract}

 \subjclass[2000]{Primary 11R06}

\keywords{Mahler measure, Lehmer problem, house of algebraic integer
}

\maketitle

\pagestyle{myheadings} \markboth{\textit{Mahler measure of
algebraic numbers}}{\textit{Chris Smyth}}

\section{Introduction}

Let
$P(x)=a_0z^d+\dots+a_d=a_0\prod_{i=1}^d(z-\alpha_i)$ be a
nonconstant polynomial with (at first) complex coefficients. Then, following
Mahler \cite{Mahl62} its {\it Mahler measure} is defined to be
\be\label{E-0}
 M(P):=\exp\left(\int_0^1 \log|P(e^{2\pi it})|dt\right),
 \ee
 the geometric mean of $|P(z)|$ for $z$ on the
unit circle. However $M(P)$ had appeared earlier in a paper of
Lehmer \cite{Lehm33}, in an alternative form
\begin{equation}\label{E-1}
M(P)=|a_0|\prod_{|\alpha_i|\ge 1} |\alpha_i|.
\ee
 The equivalence of
the two definitions follows immediately from Jensen's formula
\cite{Jens1899}
\[
\int_0^1 \log|e^{2\pi i t}-\alpha|dt=\log_+|\alpha|.
\]
Here $\log_+x$ denotes $\max(0,\log x)$. If $|a_0|\ge 1$, then clearly $M(P)\ge 1$. This is the case when $P$ has integer coefficients; we assume
 henceforth that $P$ is of this form. Then, from a result of
Kronecker \cite{Kron1857}, $M(P)=1$ occurs only if $\pm P$ is a
power of $z$ times a cyclotomic polynomial.

In \cite{Mahl62} Mahler called $M(P)$ the  measure of the polynomial
$P$, apparently to distinguish it from its (na\"\i ve) height. This
was first referred to as Mahler's measure by Waldschmidt
\cite[p.21]{Wald79} in 1979 (`mesure de Mahler'), and soon
afterwards by Boyd \cite{Boyd81} and Durand \cite{Dura81},  in the
sense of ``the function that Mahler called `measure'\,", rather than as
a name. But it soon {\it became} a name. In 1983 Louboutin \cite{Loub83}
used the term to apply to an algebraic number. We shall follow this
convention too
--- $M(\al)$ for an algebraic number
 $\al$ will mean the Mahler measure of the minimal polynomial $P_\al$ of $\al$,
  with $d$ the degree of $\al$, having conjugates $\al=\al_1,\al_2,\dots,\al_d$. The Mahler measure is
 actually a height
 function on polynomials with integer coefficients, as there are only a finite number of such polynomials of bounded degree and bounded Mahler measure.
 Indeed, in the  MR review of \cite{Loub83}, it is called the Mahler height;
 but `Mahler measure' has stuck.

 For the Mahler measure in the form $M(\al)$, there is a third
 representation to add to (\ref{E-0}) and (\ref{E-1}). We consider
 a complete set of inequivalent valuations $|.|_\nu$  of the field
 $\Q(\al)$,   normalised so that, for $\nu|p$,
$|.|_\nu=|.|_p$ on $\Q_p$. Here $\Q_p$ is the field of $p$-adic numbers, with the usual valuation $|.|_p$.
 Then for $a_0$ as in (\ref{E-1}),
\be\label{E-4}
|a_0|=\prod_{p<\infty}|a_0|^{-1}_p=\prod_{p<\infty}\prod_{\nu|p}\max(1,|\al|_\nu^{d_\nu}),
\ee coming from the product formula, and from considering the Newton polygons of the irreducible
factors (of degree $d_\nu$)  of $P_\al$ over $\Q_p$ (see e.g.
\cite[p. 73]{Weis63}).

 Then \cite[pp. 74--79]{Wald00}, \cite{BG06} from (\ref{E-1}) and (\ref{E-4})
 \be\label{E-2}
 M(\al)=\prod_{\text{all }\nu} \max(1,|\al|^{d_\nu}_\nu),
 \ee
and so also \be\label{E-5}
 h(\al):= \frac{\log M}{d}=\sum_{\text{all
 }\nu}\log_+|\al|_\nu^{d_\nu/d}.
 \ee
Here $h(\al)$ is called the {\it Weil}, or {\it absolute} height of $\al$.

\section{Lehmer's problem}
 While Mahler presumably had applications of his measure to transcendence in
mind, Lehmer's interest was in finding large primes. He sought them
amongst the Pierce numbers $\prod_{i=1}^d(1\pm\alpha_i^m)$, where
the $\al_i$ are the roots of a monic polynomial $P$ having integer
coefficients. Lehmer showed that for $P$ with no roots on the unit
circle these numbers grew with $m$ like $M(P)^m$. Pierce
\cite{Pier16} had earlier considered the factorization of these
numbers.  Lehmer posed the problem of whether, among those monic
integer polynomials with $M(P)>1$, polynomials could be chosen with
 $M(P)$ arbitrarily close to $1$. This has become known as `Lehmer's
problem', or `Lehmer's conjecture', the `conjecture' being that they
could not, although Lehmer did not in fact make this conjecture.\footnote{`Lehmer's conjecture' is also used to refer to a
conjecture on the non-vanishing of Ramanujan's $\tau$-function. But
I do not know that Lehmer actually made that conjecture either: in
\cite[p. 429]{Lehm47} he wrote ``$\dots$ and it is natural to ask
whether $\tau(n)=0$ for any $n>0$."}  The smallest value of  $M(P)>1$
he could find was
\[\label{E-1.176}
M(L)=1.176280818\dots,
\]
where $L(z)=z^{10}+z^9-z^7-z^6-z^5-z^4-z^3+z+1$ is now called `Lehmer's polynomial'. To this day no-one
 has found a smaller value of $M(P)>1$ for $P(z)\in \Z[z]$.

Lehmer's problem is central to this survey. We concentrate on
results for $M(P)$ with $P$ having integer coefficients. We do not
attempt to survey results for $M(P)$ for $P$  a polynomial in
several variables. For this we refer the reader to
 \cite{Bert07}, \cite{Boyd81}, \cite{Vill99}, \cite{Boyd98}, \cite{BV02}, \cite{Boyd02}, \cite[Chapter 3]{EW99}. However,
  the one-variable case should not really be separated from the general case, because of the fact that for every
  $P$  with integer coefficients, irreducible and in genuinely more than one variable (i.e., its Newton polytope
  is not one-dimensional) $M(P)$ is known \cite[Theorem 1]{Boyd81} to be the limit of
 $\{M(P_n)\}$ for some sequence $\{P_n\}$ of one-variable integer polynomials. This is part of a far-reaching
 conjecture of Boyd \cite{Boyd81} to the effect that the set of all $M(P)$ for $P$ an integer polynomial
 in any number of variables is a closed subset of the real line.

Our survey of results related to Lehmer's problem falls into three
categories. We report lower bounds, or sometimes exact infima, for
$M(P)$ as $P$ ranges over certain sets  of integer polynomials.
Depending on this set, such lower bounds can either tend to $1$ as
the degree $d$ of $P$ tends to infinity (Section \ref{S-small}),
be constant and greater than $1$ (Section \ref{S-medium}), or
increase exponentially with $d$ (Section \ref{S-big}). We also
report on computational work on the problem (Section
\ref{S-compute}).

In Sections \ref{S-house0} and \ref{S-house1} we discuss the closely-related function $\ho{\al}$ and the Schinzel-Zassenhaus conjecture. In Section \ref{S-disc} connections between Mahler measure and the discriminant are covered. In Section \ref{S-algnum} the known properties of $M(\al)$ as an
 algebraic number are outlined. Section \ref{S-count} is concerned with counting integer polynomials of given Mahler measure, while in Section \ref{S-dynam} a dynamical systems version of Lehmer's problem is presented. In Section \ref{S-variants}  variants of Mahler measure are discussed, and finally in Section \ref{S-applications} some applications of Mahler measure are given.

\section{The house $\ho{\al}$ of $\al$ and the conjecture of Schinzel and Zassenhaus}\label{S-house0}
Related to the Mahler measure of an algebraic integer $\al$ is $\ho{\al}$\,, the {\it house} of $\al$, defined as
the maximum modulus of its conjugates (including $\al$ itself). For $\al$ with $r>0$ roots of modulus greater
than $1$ we have the obvious inequality
\be\label{E-3}
M(\al)^{1/d}\le M(\al)^{1/r}\le \ho{\al}\le M(\al)
\ee
(see e.g. \cite{Boyd85}). If $\al$ is in fact a unit (which is certainly the case if $M(\al)<2$) then $M(\al)=M(\al^{-1})$
so that
\[
M(\al)\le \left(\max(\ho{\al}\,,\ho{{1/\al}}\,)\right)^{d/2}.
\]

In 1965 Schinzel and Zassenhaus \cite{SZ65} proved that if $\alpha\neq 0$  an algebraic integer that is not a root of unity
 and if $2s$
  of its conjugates are nonreal, then
  \be\label{E-SZ1}
  \ho{\alpha}>1+4^{-s-2}.
  \ee
  This was the first unconditional result towards solving Lehmer's problem, since by (\ref{E-3}) it implies the same lower
  bound for $M(\al)$ for such $\al$.
      They conjectured, however, that a much stronger
  bound should hold: that  under these conditions in fact
 \be\label{E-SZ}
 \ho{\alpha}\ge 1+c/d
 \ee
 for some absolute constant $c>0$.
  Its truth is implied by a positive answer to Lehmer's
   `conjecture'.  Indeed, because $\ho{\al}\ge M(\al)^{1/d}$ where $d=\deg\al$, we have
   \be\label{E-SZ2}
   \ho{\al}\ge 1+\frac{\log M(\al)}{d}=1+h(\al),
   \ee
   so that if $M(\al)\ge c_0>1$ then $\ho{\al}> 1+\frac{\log(c_0)}{d}$.

Likewise, from this inequality any results  in the direction of
   solving Lehmer's problem will have a corresponding
   `Schinzel-Zassenhaus conjecture' version. In particular, this applies to the results of Section \ref{SS-nonrec} below, including that of Breusch. His inequality appears to be the first, albeit conditional, result in the direction
   of the Schinzel-Zassenhaus conjecture or the Lehmer problem.

\section{Unconditional lower bounds for $M(\al)$ that tend to $1$\\ as $d\to\infty$}\label{S-small}
\subsection{The bounds of Blanksby and Montgomery, and Stewart.}

 The lower bound for $M(\al)$ coming from (\ref{E-SZ1}) was dramatically improved in 1971 by Blanksby and Montgomery \cite{BM71}, who showed, again for
  $\al$ of degree $d>1$ and not a root of unity, that
 \[
 M(\al)>1+\frac{1}{52d\log(6d)}.
 \]
 Their methods were based on Fourier series in several variables, making use of the nonnegativity of Fej\' er's
 kernel
\bes
\tfrac{1}{2}+\sum_{k=1}^K\left(1-\tfrac{k}{K+1}\right)\cos(kx)=\tfrac{1}{2(K+1)}
\left(\sum_{j=0}^Ke^{ix\left(\frac{K}{2}-j\right)}\right)^2. \ees
They also employed a neat geometric lemma
 for bounding the modulus of complex numbers near the unit circle: if $0<\rho\le 1$ and $\rho\le|z|\le \rho^{-1}$
  then
 \be
 |z-1|\le\rho^{-1}\left|\rho\tfrac{z}{|z|}-1\right|.
 \ee

 In 1978 Stewart \cite{Stew78a} caused some surprise by obtaining a lower
 bound of the same strength $1+\frac{C}{d\log d}$ by the use of a completely different argument. He
 based his proof on the construction of an auxiliary function of the type used in transcendence proofs.

  In such arguments it is of course necessary to make use of some arithmetic information, because of the fact that
  the polynomials one is dealing with, here the minimal polynomials of
  algebraic integers, are monic, have integer coefficients, and no root is a root of unity. In the three
  proofs of the results given above, this is done by making use of the fact that, for $\al$ not a root of unity,
  the Pierce numbers $\prod_{i=1}^d(1-\al_i^m)$ are then nonzero integers for all $m\in\N$. Hence they are at least $1$ in modulus.

\subsection{Dobrowolski's lower bound.}  In 1979 a breakthrough was achieved by  Dobrowolski, who, like Stewart, used an argument based on an
  auxiliary function to get a  lower bound
  for $M(\al)$. However, he also employed more powerful arithmetic information: the fact that for any prime $p$
  the resultant  of the minimal polynomials of $\al$ and of $\al^p$ is an integer
  multiple of $p^d$. Since this
  can be shown to be nonzero for $\al$ not a root of unity, it is at least $p^d$ in modulus. Dobrowolski
  \cite{Dobr79} was
  able to apply this fact to obtain for $d\ge 2$ the much improved lower bound
  \be\label{E-D}
  M(\al)> 1 +\frac{1}{1200} \left(\frac{\log\log d}{\log d}\right)^3.
\ee
  He also has an asymptotic version of his result, where  the constant $1/1200$ can be increased to $1-\eps$ for $\al$ of degree $d\ge d_0(\eps)$.
  Improvements in the constant in Dobrowolski's Theorem have been made since that time. Cantor and Straus \cite{CS82} proved the asymptotic version of his result with the larger constant $2-\eps$, by a different method: the auxiliary function was  replaced
  by the use of generalised Vandermonde determinants. See also \cite{Raus85} for a similar argument
  (plus some early references to these determinants). As with Dobrowolski's argument, the large
  size of the resultant of $\al$ and $\al^p$ was an essential ingredient.
   Louboutin \cite{Loub83} improved the constant further, to $9/4-\eps$, using the Cantor-Straus method. A different proof of Louboutin's result was given by Meyer \cite{Meye88}.  Later Voutier \cite{Vout96},
   by a very careful argument based on Cantor-Straus,  has obtained the constant $1/4$ valid for all $\al$ of
   degree $d\ge 2$.
    However, no-one has been able to improve the
    dependence on the degree $d$ in (\ref{E-D}), so that Lehmer's problem remains unsolved!

\subsection{Generalisations of Dobrowolski's Theorem.}
   Amoroso and David \cite{AD98,AD99} have generalised Dobrowolski's
   result in the following way. Let $\al_1,\dots,\al_n$ be $n$
   multiplicatively independent algebraic numbers in a number field
   of degree $d$. Then for some constant $c(n)$ depending only on $n$
   \be
   h(\al_1)\dots h(\al_n)\ge \frac{1}{d \log(3d)^{c(n)}}.
   \ee
Matveev \cite{Matv99} also has a result of this type, but using
instead the modified Weil height
$h_*(\al):=\max(h(\al),d^{-1}|\log\al|)$.

Amoroso and Zannier \cite{AZ} have given a version of Dobrowolski's result for $\al$, not $0$ or a root of unity, of degree $D$ over an finite abelian extension of a number field. Then
\be
h(\al)\ge \frac{c}{D}\left(\frac{\log\log{5D}}{\log{2D}}\right)^{13},
\ee
where the constant $c$ depends only on the number field, not on its abelian extension. Amoroso and Delsinne \cite {ADe} have recently improved this result,
for instance essentially reducing the exponent $13$ to $4$.

 Analogues of Dobrowolski's Theorem have been proved for elliptic curves by
 Anderson and Masser \cite{AM80}, Hindry and Silverman \cite{HS90}, Laurent \cite{Laur83} and Masser \cite{Mass89}.
  In particular  Masser proved that for an elliptic curve $E$  defined over a number field $K$ and a nontorsion
 point $P$ defined over a degree $\le d$ extension
 $F$ of $K$ that the canonical height $\hat h(P$) satisfies
 $$
 \hat h(P)\ge \frac{C}{d^{3}(\log d)^2}.
 $$
  Here $C$ depends only on $E$ and $K$. When $E$ has non-integral $j$-invariant
 Hindry and Silverman improved this bound to  $\hat h(P)\ge \frac{C}{d^{2}(\log d)^2}$. In the case where $E$ has complex multiplication, however,
  Laurent obtained the stronger bound
  $$
  \hat h(P)\ge \frac{C}{d}(\log\log d/\log d)^3.
  $$
 This is completely analogous to the formulation of Dobrowolski's result (\ref{E-D}) in terms of the Weil
  height $h(\al)=\log M(\al)/d$.

 \section{Restricted results of Lehmer strength: $M(\al)>c>1$.}\label{S-medium}

\subsection{Results for nonreciprocal algebraic numbers and polynomials}\label{SS-nonrec}
 Recall that a polynomial $P(z)$ of degree $d$ is said to be {\it reciprocal} if it satisfies $z^dP(1/z)= \pm P(z)$. (With the negative
 sign,
  clearly $P(z)$ is divisible by $z-1$.) Furthermore an algebraic number $\al$ is {\it reciprocal} if it is conjugate to
  $\al^{-1}$  (as then $P_\al$ is a reciprocal polynomial). One might at first think that it should be possible to prove
  stronger   results   on
  Lehmer's problem if we restrict our attention to reciprocal polynomials. However, this is far from being the
  case:
  reciprocal polynomials seem to be the most difficult to work with, perhaps because cyclotomic polynomials  are
  reciprocal; we can prove stronger results on Lehmer's problem if we restrict our attention to nonreciprocal
  polynomials!

 The first result in this direction was due to Breusch \cite{Breu51}. Strangely, this paper was unknown to number
  theorists until it was recently
 unearthed by Narkiewicz.
  Breusch proved that for $\al$ a nonreciprocal algebraic integer 
  \begin{equation}\label{E-Br}
  M(\al)\ge M(z^3-z^2-\tfrac{1}{4})=1.1796\dots\quad.
  \end{equation}
 Breusch's argument is based on the study of the resultant of $\al$ and $\al^{-1}$, for $\al$ a
  root of
   $P$. On the one hand, this resultant must be at least $1$ in modulus. But, on the other hand, this is not possible
   if $M(P)$ is too
   close
   to $1$, because  then all the distances $|\al_i-\overline{\al_i^{-1}}|$ are too small. (Note that $\al_i= \overline{\al_i^{-1}}$ implies that $P$ is
    reciprocal.)

 In 1971   Smyth \cite{Smyt71} independently improved the constant in (\ref{E-Br}), showing for $\al$ a nonreciprocal algebraic integer
 \be\label{E-ta}
 M(\al)\ge M(z^3-z-1)=\ta_0=1.3247\dots,
 \ee
  the real root of $z^3-z-1=0$. This constant is best possible here,
  $z^3-z-1$ being nonreciprocal. Equality $M(\al)=\ta_0$
  occurs only for $\al$ conjugate to $(\pm \ta_0)^{\pm 1/k}$
 for $k$ some positive integer.\footnote{As Boyd \cite{Boyd86b} pointed out, however,
 this does not
  preclude the possibility of equality for some {\it
  reciprocal} $\al$. But it was proved by Dixon and Dubickas
  \cite[Cor. 14]{DD04} that this could not happen.} Otherwise in fact $M(\al)>\ta_0+10^{-4}$ (\cite{Smyt72}), so that $\ta_0$
 is an isolated point in the spectrum of Mahler measures of nonreciprocal
 algebraic integers. The lower bound $10^{-4}$ for this gap in the spectrum was increased to $0.000260\dots$ by Dixon and Dubickas \cite[Th. 15]{DD04}. It would be interesting to know more about this spectrum. All of its known small
  points come from trinomials, or their irreducible factors:

 $1.324717959\dots =M(z^3-z-1)=M(\frac{z^5-z^4-1}{z^2-z+1})$;

 $1.349716105\dots=M(z^5-z^4+z^2-z+1)=M(\frac{z^7+z^2+1}{z^2+z+1})$;

 $1.359914149\dots=M(z^6-z^5+z^3-z^2+1)=M(\frac{z^8+z+1}{z^2+z+1})$;

 $1.364199545\dots=M(z^5-z^2+1)$;

 $1.367854634\dots=M(z^9-z^8+z^6-z^5+z^3-z+1)=M(\frac{z^{11}+z^4+1}{z^2+z+1})$.

 The smallest known limit point of nonreciprocal measures is
 $$\lim_{n\to\infty}M(z^n+z+1)=1.38135\dots$$
 (\cite{Boyd78b}).
 The spectrum clearly contains the set of all Pisot numbers, except perhaps  the reciprocal ones.
 But in fact it does contain those too,
 a result due to Boyd \cite[Proposition 2]{Boyd86b}. There are
 however smaller limit points of reciprocal measures (see \cite{Boyd81}, \cite{BM05}
 ).

   The method of proof of (\ref{E-ta}) was based on the Maclaurin expansion of the rational function $F(z)=P(0)P(z)/z^dP(1/z)$,
    which has integer coefficients and is nonconstant for $P$ nonreciprocal. This idea had been used in 1944
    by Salem \cite{Sale44} in his proof that the set of Pisot numbers is closed, and in the same year by
    Siegel \cite{Sieg44} in his proof that $\ta_0$ is the smallest Pisot number. One can write $F(z)$ as a
    quotient $f(z)/g(z)$ where $f$ and $g$ are both holomorphic and bounded above by $1$ in modulus in the
    disc $|z|<1$. Furthermore, $f(0)=g(0)=M(P)^{-1}$. These functions were first studied by
  Schur \cite{Schu17}, who completely specified the conditions on the coefficients of a power series $\sum_{n=0}^\infty c_nz^n$ for it to belong
   to this class. Then study of functions of this type, combined
   with the fact that the series of their quotient has integer coefficients, enables one to get the required
    lower
   bound for $M(P)$. To prove that $\ta_0$ is an isolated point of
   the nonreciprocal spectrum, it was necessary to consider the quotient
   $F(z)/F_1(z)$, where  $F_1(z)=P_1(0)P_1(z)/z^dP_1(1/z)$. Here
   $P_1$ is chosen as the minimal polynomial of some $(\pm \ta_0)^{\pm
   1/k}$ so that, if $F(z)=1+a_kz^k+\dots$, where $a_k\ne 0$ then also $F_1(z)\equiv 1+a_kz^k\pmod {z^{k+1}}$.
      Thus  this quotient, assumed nonconstant, had a first nonzero term of
   higher order, enabling one to show that $M(P)>\ta_0+10^{-4}$.

\subsection{Nonreciprocal case: generalizations}
 Soon afterwards Schinzel \cite{Schi73} and then Bazylewicz \cite{Bazy76} generalised Smyth's result to polynomials
  over
   Kroneckerian  fields.  (These are fields that are either totally real extensions of the rationals, or  totally nonreal quadratic extensions of such  fields.) For a further generalisation to polynomials in several variables see
   \cite[Theorem 70]{Schi00}. In these generalisations the optimal
   constant is obtained. If the field does not contain a primitive
   cube root of unity $\om_3$ then the best constant is again
   $\ta_0$, while if it does contain $\om_3$ then the best constant
   is the maximum modulus of the roots  $\ta$ of
   $\ta^2-\om_3\ta-1= 0$.

 Generalisations to algebraic numbers were proved by Notari \cite{Nota78} and Lloyd-Smith \cite{Lloy85}.
 See also Skoruppa's Heights notes \cite{Skor99} and  Schinzel
 \cite{Schi00}.

 \subsection{The case where {\bf $\Q(\al)/\Q$} is Galois}
 In 1999 Amoroso and David \cite{AD99}, as a Corollary of a far more general result concerning heights of
  points on subvarieties of $\G_{\text{m}}^n$, solved Lehmer's
  problem for $\Q(\al)/\Q$  a Galois extension: they  proved that
 there is a constant $c>1$ such that if $\al$ is not zero or a root of unity and $\Q(\al)$ is Galois of
  degree $d$ then $M(\al)\ge c$.

 \subsection{Other restricted results of Lehmer strength} Mignotte \cite[Cor. 2]{Mign78} proved that if $\al$
 is an algebraic number of degree $d$ such that there is a
  prime less than $d\log d$ that is unramified in the field $\Q(\al)$ then $M(\al)\ge 1.2$.

  Mignotte \cite[Prop. 5]{Mign78} gave a very short proof, based on an idea of Dobrowolski, of the fact that for
  an irreducible noncyclotomic polynomial $P$ of length $L=||P||_1$ that
  $M(P)\ge 2^{1/2L}$. For a similar result (where $2^{1/2L}$ is replaced
  by $1+1/(6L)$), see Stewart \cite{Stew78b}.

 In 2004 P. Borwein, Mossinghoff and Hare \cite{BHM04}  generalised the argument in \cite{Smyt71} to nonreciprocal
  polynomials $P$ all of whose coefficients are odd, proving that in this case
  $$M(P)\ge M(z^2-z-1)=\phi.$$  Here $\phi=(1+\sqrt{5})/{2}$. This lower bound is clearly
  best possible. Recently  Borwein, Dobrowolski and
  Mossinghoff have been able to drop the requirement of nonreciprocality: they proved
  in \cite{BDM07} that for a noncyclotomic irreducible polynomial with all odd
  coefficients then \be M(P)\ge 5^{1/4}=1.495348\dots\quad. \ee In the other direction,
  in a search \cite{BHM04}  of polynomials up to degree $72$ with coefficients $\pm 1$ and no cyclotomic factor
  the smallest  Mahler measure found was $M(z^6+z^5-z^4-z^3-z^2+z+1)=1.556030\dots$\,\,.

Dobrowolski, Lawton  and Schinzel \cite{DLS83} first
gave  a bound for the Mahler measure of an
 noncyclotomic integer polynomial $P$ in terms of the number $k$ of its nonzero coefficients: \be M(P)\ge 1+\frac{1}{\exp_{k+1}2k^2}. \ee Here
$\exp_{k+1}$ is the $(k+1)$-fold exponential.
 This was later improved by Dobrowolski \cite{Dobr91} to $1+\frac{1}{13911}\exp(-2.27k^k)$, and lately \cite{Dobr06} to
\be M(P)\ge 1+\frac{1}{\exp(a 3^{\lfloor(k-2)/4\rfloor}k^2\log k)}, \ee
where $a<0.785$.
Furthermore, in the same paper he proves that if $P$ has no  cyclotomic factors then
\be M(P)\ge 1+\frac{0.31}{k!}. \ee

With the additional restriction that $P$ is irreducible, Dobrowolski \cite{Dobr80} gave the lower bound
\be M(P)\ge
1+\frac{\log(2e)}{2e}(k+1)^{-k}.\ee
 In \cite{Dobr06} he strengthened this to
  \be M(P)\ge 1+\frac{0.17}{2^mm!}, \ee
 where $m=\lceil k/2\rceil$.

  Recently Dobrowolski \cite{Dobr07} has proved that for an integer symmetric $n\times n$ matrix $A$ with
  characteristic polynomial $\chi_A(x)$, the reciprocal polynomial
   $z^n\chi_A(z+1/z)$ is either cyclotomic or has Mahler measure at least $1.043$.
        The Mahler measure of $A$ can then be
    defined
   to be the Mahler measure of this polynomial.  McKee and Smyth \cite{MS07} have just improved the lower bound
   in Dobrowolski's result to the best possible value $\tau_0=1.176\dots$ coming from Lehmer's polynomial.
   The adjacency matrix of the graph below is an example of a matrix where this value is attained.

   The Mahler measure of a graph, defined as the Mahler measure of its adjacency matrix, has been studied by
    McKee and Smyth \cite{MS05}. They showed that  its Mahler measure was either $1$ or    at
     least $\tau_0$, the Mahler measure of the graph  
\leavevmode
\hbox{%
\epsfxsize0.8in \epsffile{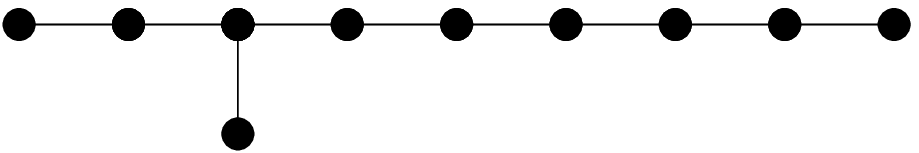}}.
   They further found all numbers in the interval $[1,\phi]$ that were Mahler measures of graphs.
   All but one of these numbers is a Salem number.

  \section{Restricted results where $M(\al)>C^d$.}\label{S-big}

\subsection{Totally real {\bf $\al$}}
  Suppose that $\al$ is a totally real algebraic integer of degree $d$, $\al\ne 0$ or $\pm 1$. Then Schinzel \cite{Schi73} proved that
  \be\label{E-S}
   M(\al)\ge \phi^{d/2}.
   \ee
 A one-page proof of this result was later provided by H\" ohn and Skoruppa
 \cite{HS93}.  The result also holds for any nonzero algebraic number $\al$ in a Kroneckerian
  field, provided $|\al|\ne 1$. Amoroso
and Dvornicich \cite[p. 261]{ADv00} gave the interesting example of $\al=\tfrac{1}{2}\sqrt{3+\sqrt{-7}}$, not an
 algebraic integer, where $|\al|=1$,
 $\Q(\al)$ is Kroneckerian, but $M(\al)=2<\phi^2$.

 Smyth \cite{Smyt80-1} studied the spectrum of values $M(\al)^{1/d}$ in $(1,\infty)$.
 He showed that this spectrum was discrete at first, and found its smallest four points. The method used is
 semi-infinite
 linear programming (continuous real variables and a finite number
 of constraints), combined with resultant information. One takes a list of judiciously chosen polynomials $P_i(x)$,
 and then finds the largest $c$ such that for some $c_i\ge 0$
 \be
 \log_+|x|\ge c -\sum_ic_i\log|P_i(x)|
 \ee
 for all real $x$. Then, averaging this inequality over the conjugates of $\al$,
 one gets that $M(\al)\ge e^c$, unless some $P_i(\al)=0$.

 Two further isolated
  points were later found by Flammang \cite{Flam96}, giving the six points comprising the whole of the spectrum
  in $(1, 1.3117)$. On the
  other hand Smyth also showed that this spectrum was dense in $(\ell,\infty)$, where $\ell=1.31427\dots$\,\,. The
  number
   $\ell$ is $\lim_{n\to\infty}M(\al_n)$, where $\beta_0=1$ and $\beta_n$, of degree $2^n$, is defined by $\beta_n-\beta_n^{-1}=\beta_{n-1} (n\ge 1)$. The limiting
    distribution of the conjugates of $\beta_n$ was studied in detail by Davie and Smyth \cite{DS89}. It is
    highly irregular: indeed, the Hausdorff dimension of the associated probability measure is
     $0.800611138269168784\dots$\,\,. It is the invariant measure of the map $\C\to\C$ taking  $t\mapsto t-1/t$,
      whose Julia set (and thus the support of the measure) is $\R$.

  Bertin \cite{Bert99} pointed out that from a result of Matveev
  (\ref{E-S}) could be strengthened when $\al$ was a nonunit.

\subsection{Langevin's Theorem}   In 1988 Langevin \cite{Lang88} proved the following general result, which
 included
Schinzel's result
   (\ref{E-S}) as a special case (though not with the explicit and best constant given by Schinzel). Suppose that
   $V$ is an
    open subset of $\mathbb C$ that has nonempty intersection with the unit circle $|z|=1$, and is stable under
    complex conjugation. Then there is a constant $C(V)>1$ such that for every irreducible monic integer
    polynomial $P$ of degree $d$ having all its roots outside $V$ one has $M(P)>C(V)^d$. The proof is based on the
     beautiful
     result of Kakeya to the effect that, for a compact subset of $\C$ stable under complex conjugation and of transfinite
     diameter
     less than $1$ there is a nonzero polynomial with integer coefficients whose maximum modulus on this set is
    less
     than $1$. (Kakeya's result is applied to the unit disc with $V$ removed.)  For Schinzel's result
    take $V=\C\setminus \R$, $C(\R)=\phi^{1/2}$, where the value of $C(\R)$ given here is best possible.
    It is of course of interest to find such best possible constants for other sets $V$.

   Stimulated by Langevin's Theorem, Rhin and Smyth \cite{RS95} studied the case where the subset of $\C$ was the
    sector $V_\ta=\{z\in\C: |\arg z|> \ta\}$.
     They found a value $C(V_\ta)>1$ for $0\le\ta\le 2\pi/3$, including $9$ subintervals of this range for which
     the constants found were best possible. In particular, the best constant $C(V_{\pi/2})$ was evaluated. This implied
      that for $P(z)$ irreducible, of degree $d$,  having all its roots with positive real part and not equal to
      $z-1$ or $z^2-z+1$ we have
      \be
     M(P)^{1/d}\ge M(z^6-2z^5+4z^4-5z^3+4z^2-2z+1)^{1/6}=1.12933793\dots,
     \ee
      all roots of $z^6-2z^5+4z^4-5z^3+4z^2-2z+1$ having positive real part. Curiously, for
     some  root $\al$ of this polynomial, $\al+1/\al=\ta_0^2$, where as above $\ta_0$ is the smallest Pisot
      number.

     Recently Rhin and Wu \cite{RW05} extended these results, so that there are now $13$ known subintervals of
      $[0,\pi]$ where the best constant $C(V_\ta)$ is known.
It is of interest to see what happens as $\ta$ tends to $\pi$; maybe one could obtain a bound connected to
 Lehmer's original problem. Mignotte \cite{Mign89}
has looked at this, and has shown that for $\ta=\pi-\eps$ the smallest limit point of the set $M(P)^{1/d}$
for $P$ having all its roots outside $V_\ta$ is at least $1+c\eps^3$ for some positive constant $c$.

Dubickas and Smyth \cite{DS01a} applied Langevin's Theorem to the annulus
$$V(R^{-\gamma},R)=\{z\in\C\mid R^{-\gamma}<|z|<R\},$$ where $R>1$ and $\gamma>0$, proving that the best constant
$C(V(R^{-\gamma},R))$ is $R^{\gamma/(1+\gamma)}$.

 \subsection{Abelian number fields}  In 2000 Amoroso and Dvornicich \cite{ADv00} showed that when $\al$ is a nonzero algebraic number, not a root of unity,
 and $\Q(\al)$
  is an abelian extension of $\Q$ then $M(\al)\ge 5^{d/12}$. They also give an example with $M(\al)=7^{d/12}$. It would be interesting to find the best constant $c>1$ such that $M(\al)\ge c^d$ for these numbers.
Baker and Silverman \cite{Bake03}, \cite{Silv04}, \cite{BaSi04}
generalised this lower bound first to elliptic curves, and then to
abelian varieties of arbitrary dimension.

 \subsection{Totally {\bf $p$}-adic fields} Bombieri and Zannier \cite{BZ01} proved an analogue of Schinzel's result (\ref{E-S})
 for `totally $p$-adic' numbers:
   that is, for algebraic numbers $\al$ of degree $d$ all of whose conjugates lie in $\Q_p$. They showed that then  $M(\al)\ge c_p^d$, for some constant $c_p>1$.

\subsection{The heights of Zagier and Zhang and generalisations.}
Zagier \cite{Zagi93} gave a result that can be formulated as
proving that the Mahler measure of any irreducible nonconstant
polynomial in $\Z[(x(x-1)]$ has Mahler
measure at least $\phi^{d/2}$, apart from  $\pm(x(x-1)+1)$. Doche \cite{Doch01a,Doch01b} studied
the spectrum resulting from the measures of such polynomials, giving
a gap to the right of the smallest point $\phi^{1/2}$, and finding a
short interval where the smallest limit point lies. He used the
semi-infinite linear programming method outlined above. For this
problem, however, finding the second point of the spectrum seems to
be difficult. Zagier's work was motivated by a far-reaching result
of Zhang \cite{Zhan92} (see also \cite[p. 103]{Wald00}) for curves
on a linear torus. He proved that for all such curves, apart from those of the type
$x^iy^j=\om$, where $i,j\in\Z$ and $\om$ is a root of unity, there is a constant $c>0$ such that the curve has only
finitely many algebraic points $(x,y)$ with $h(x)+h(y)\le c$.
Zagier's result was for the curve $x+y=1$.

Following on from Zhang, there have been recent deep and diverse
generalisations in the area of small points on subvarieties of
$\G_{\text{m}}^n$. In particular see Bombieri and Zannier
\cite{BZ95}, Schmidt \cite{Schm96} and Amoroso and David
 \cite{AD00,AD03,AD04,AD06}.

Rhin and Smyth \cite{RS97} generalised Zagier's result  by replacing
 polynomials in $\Z(x(x-1))$ by polynomials in $\Z[Q(x)]$, where
$Q(x)\in\Z[x]$ is not $\pm$ a power of $x$. Their proof used a very
general result of Beukers and Zagier \cite{BZ97} on heights of points
on projective hypersurfaces.
  Noticing that Zagier's result has the
  same lower bound as Schinzel's result above for totally real $\al$, Samuels \cite{Samu06} has recently
  shown that
  the same lower bound holds for a more general height
  function. His result  includes  those of both Zagier and Schinzel.
  The proof is also based on \cite{BZ97}.

\section{Lower bounds for $\ho{\al}$}\label{S-house1}

\subsection{General lower bounds} We know that any lower bound for $M(\al)$ immediately gives a corresponding lower
bound for $\ho{\al}$\,, using (\ref{E-SZ2}). For instance, from
  \cite{Vout96}  it follows that for $\al$ of degree $d>2$ and not a root of unity
 \be\label{E-V}
 \ho{\al}\ge 1+\frac{1}{4d}\left(\frac{\log\log d}{\log d}\right)^3.
 \ee
 Some lower bounds, though asymptotically weaker, are better for small degrees.
 For example Matveev \cite{Matv91} has shown that for such $\al$
 \be\label{E-M}
 \ho{\al}\ge \exp\frac{\log(d+0.5)}{d^2},
 \ee
which is better than (\ref{E-V}) for $d\le 1434$ (see \cite{RW07}).
Recently Rhin and Wu have improved (\ref{E-M}) for $d\ge 13$ to
 \be
 \ho{\al}\ge \exp\frac{3\log(d/2)}{d^2},
 \ee
which is better than (\ref{E-V}) for $d\le 6380$. See also the paper of Rhin and Wu in this volume.

 Matveev \cite{Matv91} also proves that if $\al$ is a reciprocal
  (conjugate to $\al^{-1}$) algebraic integer, not a root of unity,
  then $\ho{\al}\ge (p-1)^{1/(pm)}$, where $p$ is the least prime
  greater than $m=n/2\ge 3$.

   Indeed, Dobrowolski's first result in this area \cite{Dobr78}
   was for $\ho{\al}$ rather than $M(\al)$: he proved that
   \bes
   \ho{\al}> 1+\frac{\log d}{6d^2}.
   \ees
   His argument is a beautifully simple one,
    based on the use of the power sums
   $s_k=\sum_{i=1}^d \al_i^k$, the Newton identities, and the arithmetic fact that, for any prime $p$,
    $s_{kp}\equiv s_k\pmod p$.

    The strongest asymptotic result to date in the direction of the Schinzel-Zassenhaus conjecture is due
    to Dubickas \cite{Dubi93}: that
   given $\eps>0$ there is a constant $d(\eps)$ such than any
   nonzero algebraic integer $\al$ of degree $d>d(\eps)$ not a root of unity
   satisfies
   \be\label{E-64}
   \ho{\al}>
   1+\left(\frac{64}{\pi^2}-\eps\right)\left(\frac{\log\log d}{\log
   d}\right)^3\frac{1}{d}.
   \ee

Cassels \cite{Cass66} proved that if an algebraic number $\al$ of
degree $d$ has the property $\ho{\al}\le 1+\frac{1}{10d^2}$ then at least one of
the conjugates of $\al$ has modulus $1$. Although this result has
been superseded by Dobrowolski's work, Dubickas \cite{Dubi00}
applied  the inequality \be \prod_{k<j}|z_k\overline{z_j}-1|\le
n^{n/2}\left(\prod_{m=1}^n\max(1,|z_m|)\right)^{n-1} \ee for complex
numbers $z_1,\dots,z_n$, a variant of one in \cite{Cass66},
 to prove that
$$M(\al)^2\left|\prod\log|\al_i|\right|^{1/d}\ge 1/(2d)$$
 for a nonreciprocal algebraic number $\al$ of degree $d$ with conjugates $\al_i$.

\subsection{ The house {\bf$\ho{\al}$} for {\bf$\al$} nonreciprocal}\label{S-3/2}
The Schinzel-Zassenhaus conjecture (\ref{E-SZ}) restricted to
nonreciprocal polynomials follows from Breusch's result above, with
$c=\log 1.1796\dots=0.165\dots$, using (\ref{E-SZ2}). Independently
Cassels \cite{Cass66} obtained this result with $c=0.1$, improved by
Schinzel to $0.2$
 (\cite{Schi69}), and by Smyth \cite{Smyt71} to $\log \ta_0=0.2811\dots$\,\,. He also showed that $c$ could not
 exceed $\frac{3}{2}\log\ta_0=0.4217\dots$\,\,.
In 1985 Lind and Boyd (see \cite{Boyd85}), as a result of extensive computation (see Section \ref{S-compute}),
conjectured that, for  degree $d$, the extremal $\al$ are nonreciprocal and have
$\sim{\frac{2}{3}}d$ roots outside the unit circle. What a contrast with Mahler measure, where all small
$M(\al)$ are reciprocal!
 This would imply that the best constant $c$ is $\frac{3}{2}\log\ta_0$. In 1997 Dubickas \cite{Dubi97}
proved that $c>0.3096$ in this nonreciprocal case.

 \subsection{The house of totally real  {\bf $\al$}.}
 Suppose that $\al$ is a totally real algebraic integer.
 If  $\ho{\al}\le 2$ then by \cite[Theorem 2]{Kron1857} $\al$
 is of the form $\om+1/\om$, where $\om$ is a root of unity.  If for some $\delta>0$ we have $2<\ho{\al}\le 2+\delta^2/(1+\delta)$, then, on
 defining $\gamma$ by
 $\gamma+1/\gamma=\al$, we see that $\gamma$ and its conjugates are either
 real or lie on the unit circle, and $1<\ho{\gamma}\le 1+\delta$. This
 fact readily enables us to deduce a lower bound greater than $2$ for $\ho{\al}$ whenever we have a lower
 bound greater than $1$ for $\ho{\gamma}$\,. Thus from (\ref{E-SZ1}) \cite{SZ65} it
 follows that for $\al$  not of the form $2\cos \pi r$ for any $r\in \Q$
 \be\label{E-SZ1R}
\ho{\al}\ge 2+4^{-2d-3} \ee \cite{SZ65}. In a similar way
(\ref{E-64}) above implies that for such $\al$, and $d>d(\eps)$
that
 \be\label{E-64R}
   \ho{\al}>
   2+\left(\frac{4096}{\pi^4}-\eps\right)\left(\frac{\log\log d}{\log
   d}\right)^6\frac{1}{{d}^2}
   \ee
   \cite{Dubi93}.
   However Dubickas \cite{Dubi95} managed to improve this lower bound to
    \be \ho{\al}>2+3.8\frac{(\log\log d)^3}{d(\log d)^4}.
\ee He improved the constant $3.8$ to $4.6$ in \cite{Dubi97}.

\subsection{The Kronecker constant} Callahan, Newman and Sheingorn \cite{CNS77} define the {\it
Kronecker constant} of a number field $K$ to be the least $\eps>0$ such that
$\ho{\al}\ge 1+\eps$ for every algebraic integer $\al\in K$. The truth of the Schinzel-Zassenhaus conjecture
(\ref{E-SZ}) would imply that
the Kronecker constant of $K$ is at least $c/[K:\Q]$. They give  \cite[Theorem 2]{CNS77} a sufficient condition
 on $K$ for this to be the case.
They also point out,  from
considering $\al\overline{\al}-1$, that if $\al$ is a nonzero
algebraic integer not a root
    of unity in a Kroneckerian  field then $\ho{\al}\ge \sqrt{2}$ (See also \cite{Mign81}), so
 that the Kronecker constant of a Kroneckerian  field is at least $\sqrt{2}-1$.

  \section{small values of $M(\al)$ and $\ho{\al}$}\label{S-compute}

\subsection{Small values of { $ M(\al)$} }  The first recorded computations on Mahler measure were performed by Lehmer in his 1933 paper \cite{Lehm33}. He found the smallest values of $ M(\al)$ for $\al$ of degrees $2,3$ and $4$, and the smallest $M(\al)$ for $\al$
reciprocal of degrees $2,4,6$ and $8$. Lehmer records the fact that Poulet (?unpublished) ``\dots has made a similar investigation of symmetric polynomials with practically the same results". Boyd has done extensive computations, searching for `small' algebraic integers
  of various kinds. His first major published table was of  Salem
  numbers less than $1.3$
   \cite{Boyd77}, with four more found in \cite{Boyd78a}.
  Recall that these are positive reciprocal algebraic integers of degree at least $4$ having only one conjugate
  (the number itself) outside the unit circle.
These numbers give many examples of small Mahler measures, most
notably (from (\ref{E-1.176})) $M(L)=1.176\dots$ from the Lehmer
polynomial itself, which is the minimal polynomial of a Salem
number. In  later computations \cite{Boyd80}, \cite{Boyd89}, he finds all reciprocal
$\al$ with $M(\al)\le 1.3$ and degree up to $20$, and those with
$M(\al)\le 1.3$ and degree up to $32$ having coefficients in
$\{-1,0,1\}$ (`height $1$').

Mossinghoff \cite{Moss98} extended Boyd's tables from degree $20$ to
degree $24$ for $M(\al)<1.3$, and to degree $40$ for height $1$
polynomials, finding four more Salem numbers less than $1.3$.  He
also has a website  \cite{Moss} where up-to-date tables of small Salem
numbers and Mahler measures are conveniently displayed (though
unfortunately without their provenance). Flammang, Grandcolas and
Rhin \cite{FGR99} proved that Boyd's  table, with the additions by
Mossinghoff, of the $47$ known Salem numbers less than $1.3$ is
complete up to degree $40$. Recently Flammang, Rhin and Sac-\' Ep\'
ee \cite{FRS06} have extended these tables, finding all
$M(\al)<\ta_0$ for $\al$ of degree up to $36$, and all $M(\al)<1.31$
for $\al$ of degree up to $40$. This latter computation showed that
the earlier tables of Boyd and Mossinghoff for $\al$ of degree up to
$40$ with $M(\al)<1.3$ are complete.

\subsection{Small values of { $\ho{\al}\,$}} Concerning $\ho{\al}\,$, Boyd \cite{Boyd85} gives tables of
 the smallest values of $\ho{\al}$ for $\al$ of  degree $d$ up to $12$, and for $\al$ reciprocal of degree up to $16$. Further
computation has recently been done on this problem by Rhin and Wu
\cite{RW07}. They computed the smallest house of algebraic numbers of degree up to $28$. All are nonreciprocal, as predicted by Boyd's conjecture (see Section \ref{S-3/2}). Their data led the authors to conjecture that,
for a given degree, an algebraic number of that degree with minimal house was a root of a polynomial consisting of at most four monomials.

  \section{Mahler measure and the discriminant}\label{S-disc}

\subsection{}
  Mahler \cite{Mahl64b} showed that for a complex polynomial $$P(z)=a_0z^d+\dots+a_d=a_0(z-\al_1)\dots(z-\al_d)$$
  its
   discriminant $\disc(P)=a_0^{2d-2}\prod_{i<j}(\al_i-\al_j)^2$ satisfies
 \be
 |\disc(P)|\le d^dM(P)^{2d-2}.
 \ee
 From this it follows immediately that if there is an absolute constant $c>1$ such that $|\disc(P)|\ge (cd)^d$
  for all irreducible $P(z)\in\mathbb Z[z]$,
 then $M(P)\ge c^{d/(2d-2)}$, which would solve Lehmer's problem. This consequence of Mahler's inequality has
 been
 noticed in various variants by several people, including Mignotte
 \cite{Mign78} and Bertrand \cite{Bert82}.

 In 1996 Matveev \cite{Matv96} showed that in Dobrowolski's inequality, the
 degree $d\ge 2$ of $\al$ could be replaced by a much smaller (for large
 $d$) quantity $$\delta=\max(d/\disc(\al)^{1/d},\delta_0(\eps))$$ for those $\al$
 for which $\al^p$ had degree $d$ for all primes $p$. (Such $\al$ do not include any roots of unity.)
  Specifically, he obtained for given $\eps>0$
  \be
  M(\al)\ge \exp\left((2-\eps)\left(\frac{\log\log \delta}{\log
  \delta}\right)^3\right)
  \ee
for these $\al$.

 Mahler \cite{Mahl64b} also gives the lower bound
 \be
 \delta(P)>\sqrt{3}|\disc(P)|^{1/2}d^{-(d+2)/2}M(P)^{-(m-1)}
 \ee
for the minimum distance $\delta(P)=\min_{i<j}|\al_i-\al_j|$ between
the roots of $P$.

\subsection{Generalisation involving the discriminant  of Schinzel's lower
bound}\label{S-discRhin}
   Rhin \cite{Rhin04} generalised Schinzel's result (\ref{E-S}) by proving, for $\al$ a totally positive algebraic
  integer of degree $d$ at least $2$ that
  \be
  M(\al)\ge \left(\frac{\delta_1+\sqrt{\delta_1^2+4}}{2}\right)^{d/2}.
  \ee
  Here $\delta_1=|\disc(\al)|^{1/d(d-1)}$. This result apparently also follows from an earlier result of Za\"\i mi \cite{Zaim94} concerning a lower bound for a weighted product of the moduli of the conjugates of an algebraic integer --- see the Math Review of Rhin's paper.

\section{Properties of $M(\al)$ as an algebraic number}\label{S-algnum}

A {\it Perron number} is an algebraic integer with exactly one
conjugate of maximum modulus. It is clear from  (\ref{E-1}) that
$M(\al)$ is a Perron number for any algebraic integer $\al$; this
seems to have been first observed by Adler and Marcus \cite{AM79}
(see \cite{Boyd86b}). In the other direction: is the Perron number
$1+\sqrt{17}$ a Mahler measure? See Schinzel \cite{Schi04}, Dubickas
\cite{Dubi05}. Dubickas \cite{Dubi04c} proves that for any Perron
number $\beta$ some integer multiple of $\beta$ is a Mahler measure.
(These papers also contains other interesting properties of the set
of Mahler measures.) Boyd \cite{Boyd86a} proves that if
$\beta=M(\al)$ for some algebraic integer $\al$, then all conjugates
of $\beta$ other than $\beta$ itself either lie in the annulus
$\beta^{-1}<|z|<\beta$ or are equal to $\pm\beta^{-1}$.

 If $\al$ were reciprocal, it might be expected that
$M(\al)$ would be reciprocal too, while if $\al$ were nonreciprocal,
then $M(\al)$ would be nonreciprocal. However neither of these need
be the case: in \cite[Proposition 6]{Boyd86b} Boyd exhibits a family
of degree $4$ Pisot numbers that are the Mahler measures of
reciprocal algebraic integers of degree $6$, and in
\cite[Proposition 2]{Boyd86b} he notes that for $q\ge 3$ a root
$\al_q$ of the irreducible nonreciprocal polynomial
$z^4-qz^3+(q+1)z^2-2z+1$ then $M(\al_q)=\frac{1}{2}(q+\sqrt{q^2-4})$
is reciprocal.  In fact, since
$M(\frac{1}{2}(q+\sqrt{q^2-4}))=\frac{1}{2}(q+\sqrt{q^2-4})$, this
also shows that a number can be both a reciprocal and a
nonreciprocal measure.  See also \cite{Boyd87}. Dixon and Dubickas
\cite{DD04} prove that the set of all $M(\al)$ does not form a
semigroup, as for instance $\sqrt{2}+1$ and $\sqrt{3}+2$ are Mahler
measures, while their product is not. (In terms of polynomials, this
set is of course equal to the set of all $M(P)$ for $P$ irreducible.
If instead we take the set of all (reducible and irreducible)
polynomials, then, because of $M(PQ)=M(P)M(Q)$ this larger set {\it
does} form a semigroup.)

In \cite{Dubi04b} Dubickas proves that the additive group generated
by all Mahler measures is the group of all real algebraic numbers,
while the multiplicative group generated by all Mahler measures is
the group of all positive real algebraic numbers.

We know that $M(P(z))=M(P(\pm z^k))$ for either choice of sign, and
any $k\in \N$. Is this the only way that Mahler measures of
irreducible polynomials can be equal? Boyd \cite{Boyd80} gives some
illuminating examples to show that there can be other reasons that make this
happen. The examples were discovered during his computation of
 reciprocal polynomials of small Mahler measure (see Section \ref{S-compute}). For
 example, for $P_6=z^6+2z^5+2z^4+z^3+2z^2+2z+1$ and $P_8=z^8+z^7-z^6-z^5+z^4-z^3-z^2+z+1$ we have
\begin{equation*}
M(P_6)=M(P_8)=1.746793\ldots=M,
\end{equation*}
say, where both polynomials are irreducible. Boyd explains how such
examples arise. If $\al_i(i=1,\dots,8)$ are the roots of $P_8$, then
for different $i$ $M(\al_1\al_i)$ can equal $M$, $M^2$ or $M^3$. The
roots of $P_6$ are the three $\al_1\al_i$ with $M(\al_1\al_i)=M$ and
their reciprocals. Clearly $M(\al_1^2)=M^2$, while for three other
$\al_i$ the product $\al_1\al_i$ is of degree $12$ and has
$M(\al_1\al_i)=M^3$. ($P_8$ has the special property that it has
roots $\al_1, \al_2,\al_3, \al_4$ with $\al_1\al_2=\al_3\al_4\ne 1$.)

Dubickas \cite{Dubi02} gives a lower bound for the distance of an
algebraic number $\gamma$ of degree $n$ and leading coefficient $c$,
not a Mahler measure, from a Mahler measure $M(\al)$ of degree $D$:
\be |M(\al)-\gamma|>c^{-D}(2\,\ho{\gamma}\,)^{-nD}. \ee

\section{Counting polynomials with given Mahler measure}\label{S-count}
Let $\#(d,T)$ denote the number  of integer polynomials of degree
$d$ and Mahler measure at most $T$. This function has been studied
by several authors. Boyd and Montgomery \cite{BM90} give the
asymptotic formula \be c(\log
d)^{-1/2}d^{-1}\exp\left(\frac{1}{\pi}\sqrt{105\zeta(3)d}\right)(1+o(1)),
\ee where $c=\frac{1}{4\pi^2}\sqrt{105\zeta(3)e^{-\gamma}}$, for the
number $\#(d,1)$ of cyclotomic polynomials of degree $d$, as
$d\to\infty$.

 Dubickas and Konyagin \cite{DK98} obtain by simple arguments the
 lower bound $\#(d,T)>\frac{1}{2}T^{d+1}(d+1)^{-(d+1)/2}$, and
 upper bound $\#(d,T)<T^{d+1}\exp(d^2/2)$, the latter being valid
 for $d$ sufficiently large. For $T\ge \ta_0$ they derived the upper
 bound  $\#(d,T)<T^{d(1+16\log\log d/\log d)}$.
 Chern and Vaaler \cite{CV01} obtained
the asymptotic formula $V_{d+1}T^{d+1}+O_d(T^d)$ for $\#(d,T)$ for
fixed $d$, as $T\to\infty$. Here $V_{d+1}$ is an explicit constant
(the volume of a certain star body). Recently Sinclair \cite{Sinc07} has produced corresponding estimates for counting functions of reciprocal polynomials.

\section{A dynamical Lehmer's problem}\label{S-dynam}
Given a rational map $f(\al)$ of degree $d\ge 2$ defined over a number field $K$, one
can define for $\al$ in some extension field of $K$ a canonical height
$$
h_f(\al) =\lim_{n\to\infty} d^{-n}h(f^n(\al)),
$$
where $f^n$ is the $n$th iterate of $f$, and $h$ is, as before, the Weil height of $\al$. Then $h_f(\al)=0$ if and only if the iterates $f^n(\al)$ form a finite set, and an analogue of Lehmer's  problem would be to decide whether or not
$$
h_f(\al) \ge \frac{C}{\deg(\al)}
$$
for some constant $C$ depending only on $f$ and $K$. Taking $f(\al)=\al^d$ we retrieve the Weil height and the original Lehmer problem. There seem to be no good estimates, not even of polynomial decay, for any $f$ not associated to an endomorphism of an algebraic group. See \cite[Section 3.4]{Silv07} for more details.

\section{Variants of Mahler measure}\label{S-variants}
Everest and n\'{\i}  Fhlath\' uin \cite{EF96} and Everest and Pinner
\cite{EP98} (see also \cite[Chapter 6]{EW99}) have defined the {\it
elliptic Mahler measure}, based on a given elliptic curve $E=\C/L$
over $\C$, where $L=\langle \om_1,\om_2\rangle\subset \C$ is a
lattice, with $\wp_L$ its associated Weierstrass $\wp$-function. Then
for $F\in\C[z]$ the (logarithmic) elliptic Mahler measure $m_E(F)$
is defined as \be
\int_0^1\int_0^1\log|F(\wp_L(t_1\om_1+t_2\om_2))|dt_1dt_2. \ee If
$E$ is in fact defined over $\Q$ and has a rational point $Q$ with
$x$-coordinate $M/N$ then often $m_E(Nz-M)=2\hat h(Q)$, showing that
$m_E$ is connected with the canonical height on $E$.

Kurokawa \cite{Kuro04} and Oyanagi \cite{Oyan04} have defined a $q$-analogue
of Mahler measure, for a real parameter $q$. As $q\to 1$ the classical Mahler measure is recovered.

Dubickas and Smyth \cite{DS01b} defined the {\it metric Mahler
measure} ${\mathcal M}(\al)$ as the infimum of $\prod_i
M(\beta_i)$, where $\prod_i \beta_i=\al$. They used this to define
a metric on the group of nonzero algebraic numbers modulo torsion
points, the metric giving the discrete topology on this group if
and only if Lehmer's `conjecture' is true (i.e.,~
$\inf_{\al:M(\al)>1}M(\al)>1$).

Very recently Pritsker \cite{Prit08, Prit07} has studied an areal analogue of Mahler measure, defined by replacing the normalised arclength measure on the unit circle by the normalised area measure on the unit disc.

  \section{applications}\label{S-applications}

\subsection{Polynomial factorization}
  I first met Andrzej Schinzel at the ICM in Nice in 1970. There he mentioned to me an
  application of Mahler measure to irreducibility of polynomials.  (After this we had some correspondence about the work leading to \cite{Smyt71}, which was very helpful to me.) If a class of irreducible polynomials
  had Mahler measure at least $B$, then any polynomial
  of Mahler measure less than $B^2$ can have at most one factor from that class.
 For instance, a trinomial $z^d\pm z^m\pm 1$ has,  by Vicente Gon\c{c}alves'
 inequality \cite{ViGo50}, \cite{Ostr60}, \cite[Th. 9.1.1]{RS02}
 $M(P)^2+M(P)^{-2}\le || P||_2^2$, Mahler measure at most $\phi$.
 Since $\phi<\ta_0^2$, by (\ref{E-ta}) such trinomials can have at most one irreducible noncyclotomic factor.  Here $|| P||_2$ is the $2$-norm of $P$ (the square root of the sum
   of the squares of its coefficients).

  More generally   Schinzel (see \cite{Dobr79}) pointed
   out  the following consequence of (\ref{E-D}): that for any
   fixed $\eps>0$ and polynomial $P$ of degree $d$ with integer coefficients, the
   number of its noncyclotomic irreducible factors counted with multiplicities is $O(d^\eps||
   P||_2^{1-\eps})$.
  See also \cite{Schi76}, \cite{Schi83},  \cite{PV93-9}.

\subsection{Ergodic theory}
One-variable Mahler measures have applications in ergodic theory.
Consider an automorphism of the torus $\R^d/\Z^d$ defined by a
$d\times d$ integer matrix of determinant $\pm 1$, with
characteristic polynomial $P(z)$. Then the topological entropy of
this map is  $\log M(P)$ (Lind \cite{Lind74} --- see
also \cite{Boyd81}, \cite[Theorem 2.6]{EW99}).

\subsection{Transcendence and diophantine approximation}  Mahler measure, or rather the Weil height $h(\al)=\log
M(\al)/d$, plays an important technical r\^ ole
in modern transcendence theory, in particular for bounding the
coefficients of a linear form in logarithms known to be dependent.

As remarked by Waldschmidt \cite[p65]{Wald00}, the fact that this
height has three equivalent representations, coming from
(\ref{E-0}), (\ref{E-1}) and (\ref{E-2}) makes it a very versatile
height function for these applications.

If $\al_1,\dots,\al_n$ are algebraic numbers such that their logarithms are $\Q$-linearly dependent,
then it is of importance in Baker's transcendence method to get small upper estimates for the size of integers
$m_1,\dots,m_n$ needed so that $m_1\log\al_1+\dots+m_n\log\al_n=0$. Such estimates can be given using Weil
heights of the $\al_i$. See \cite[Lemma 7.19]{Wald00} and the remark after it.

Chapter 3 (`Heights of Algebraic Numbers') of \cite{Wald00} contains
 a wealth of interesting material on the Weil height and
other height functions, connections between them, and applications.
For instance, for a polynomial $f\in\Z[z]$ of degree at most $N$ for
which the algebraic number $\al$ is not a root one has \bes
|f(\al)|\ge\frac{1}{M(\al)^N|| f ||_1^{d-1}},
 \ees
  where $|| f ||_1$ is the length of
$f$, the sum of the absolute values of its coefficients, and
$d=\deg\al$ (\cite[p83]{Wald00}).

In particular,  for a rational number $p/q\ne \al$ with $q>0$, and
$f(x)=qx-p$ we obtain \be\label{E-pq}
\left|\al-\frac{p}{q}\right|\ge\frac{1}{M(\al)q(\max(|p|+q))^{d-1}}.
\ee
\subsection{Distance of $\al$ from $1$.} From (\ref{E-pq}) we immediately
get for $\al\ne 1$ \be |\al-1|\ge \frac{1}{2^{d-1}M(\al)}. \ee
Better lower bounds for $|\al-1|$ in terms of its Mahler measure
have been given by Mignotte \cite{Mign79}, Mignotte
and Waldschmidt \cite{MW94}, Bugeaud, Mignotte and
 Normandin \cite{BMN95}, Amoroso \cite{Amor96}, Dubickas \cite{Dubi95}, and  \cite{Dubi98}.
For instance Mignotte and Waldschmidt prove that
 \be
|\al-1|>\exp\{-(1+\eps)(d(\log d)(\log M(\al)))^{1/2}\} \ee for
$\eps>0$ and $\al$ of degree $d\ge d(\eps)$. Dubickas \cite{Dubi95}
improves the constant $1$ in this result to $\pi/4$, and in the
other direction \cite{Dubi98} proves that for given $\eps>0$ there
is an infinite sequences of degrees $d$ for which an $\al$ of degree
$d$ satisfies \be |\al-1|<\exp\left\{-(c-\eps)\left(\frac{d\log
M(\al)}{\log d}\right)^{1/2}\right\}. \ee Here Dubickas uses the
following simple result: if $F\in\C[z]$ has degree $t$ and $F'(1)\ne
0$ then there is a root $a$ of $F$ such that $|a-1|\le
t|F(1)/F'(1)|$.

\subsection{Shortest unit lattice vector}
Let $K$ be a number field with unit lattice of rank $r$, and $M=\min
M(\al)$, the minimum being taken over all units $\al\in K$, $\al$
not a root of unity. Kessler \cite{Kess91} showed that  then the
shortest vector $\lambda$ in the unit lattice has length
$||\lambda||_2$ at least $\sqrt{\frac{2}{r+1}}\log M$.

\subsection{Knot theory}
Mahler measure of one-variable polynomials arises in knot theory
in connection with Alexander polynomials of knots and reduced
Alexander polynomials of links --- see Silver and Williams
\cite{SW02}. Indeed, in Reidemeister's classic book on the subject
\cite{Reid32}, the polynomial $L(-z)$ appears as the Alexander
polynomial of the $(-2,3,7)$-pretzel knot. Hironaka \cite{Hiro01}
has shown that among a wide class of Alexander polynomials of
pretzel links, this one has the smallest Mahler measure.
Champanerkar and  Kofman \cite{CK05} study a sequence of Mahler
measures of Jones polynomials of hyperbolic links $L_m$ obtained
using $(-1/m)$-Dehn surgery, starting with a fixed link. They show
that it converges to the Mahler measure of a $2$-variable
polynomial. (The many more applications of Mahler measures of
several-variable polynomials to geometry and topology are outside
the scope of this survey.)

\section{Final remarks}

\subsection{Other sources on Mahler measure}

 Books covering various aspects of Mahler measure include the following:
 Bertin and Pathiaux-Delefosse \cite{BP89}, Bertin {\it et al} \cite{BDGPS92}, Bombieri and
 Gubler \cite{BG06}, Borwein \cite{Borw02},
  Schinzel \cite{Schi82}, Schinzel \cite{Schi00},  Waldschmidt \cite{Wald00}.

Survey articles   and lecture notes on Mahler measure include:
 Boyd \cite{Boyd78b}, Boyd \cite{Boyd81},  Everest \cite{Ever98},
 Hunter \cite{Hunt82}, Schinzel \cite{Schi99}, Skoruppa \cite{Skor99}, Stewart \cite{Stew78b},
 Vaaler \cite{Vaal03}, Waldschmidt \cite{Wald81}.

 \subsection{Memories of Mahler} As one of a small group of undergraduates in ANU, Canberra in
 the mid-1960s, we were encouraged to attend  graduate courses at the
 university's Institute of Advanced Studies, where Mahler had a
 research chair. I well remember his lectures on transcendence with his  blackboard copperplate handwriting,
  all the technical details being carefully spelt out.

\subsection{Acknowledgements} I thank Matt Baker, David Boyd, Art\= uras Dubickas, James McKee, Alf van der Poorten, Georges Rhin, Andrzej Schinzel,
Joe Silverman, Michel Waldschmidt, Susan Williams, Umberto Zannier and the referee
for some helpful remarks concerning an earlier draft of this survey,
which have been incorporated into the final version. This article
arose from a talk I gave in January 2006 at the Mahler Measure in
Mobile Meeting, held in Mobile, Alabama. I would like to thank the
organisers Abhijit Champanerkar,
 Eriko Hironaka, Mike Mossinghoff, Dan Silver and Susan Williams for a stimulating meeting.

\end{document}